\documentclass[12pt]{article}
\usepackage[left=2.5cm,right=2.5cm,top=2.5cm,bottom=2.5cm,a4paper]{geometry}

\usepackage{graphics}
\usepackage{graphicx}
\usepackage{epsfig}
\usepackage{amsmath,amsfonts,amssymb}
\usepackage{listings}
\usepackage{paralist}
\usepackage{natbib}
\usepackage{sectsty}
\numberwithin{equation}{section}
\date{}

\subsectionfont{\normalfont}

\newtheorem{theorem}{Theorem}[section]
\newtheorem{definition}[theorem]{Definition}
\newtheorem{lemma}[theorem]{Lemma}
\newtheorem{corollary}[theorem]{Corollary}

\newtheorem{example}[theorem]{Example}
\newtheorem{proposition}[theorem]{Proposition}

\newtheorem{open problem}[theorem]{Open problem}

\begin{document}

\title{The frobenius problem for generalized thabit numerical semigroups}

\maketitle
$\phantom{1}$\\
Kyunghwan Song$^{1,*}$\\
\textit{${}^{1}$Department of Mathematics, Korea University, 145 Anam-ro, Seongbuk-gu, Seoul, Korea\\
${}^{*}$Corresponding Author: Email: heroesof@korea.ac.kr}.\\
\\
\textbf{Abstract}\\
\\
In this paper,  we will introduce the Frobenius problem for generalized Thabit numerical semigroups, which is motivated by the Frobenius problem for Thabit numerical semigroups.

\section{Introduction}\label{sec_Introduction}
The greatest integer that does not belong to a subset $S$ of the natural integers is called the Frobenius number of $S$ and denoted by $F(S)$. In other words, the Frobenius number is the largest integer that cannot be expressed as a sum $\sum_{i=1}^n t_i a_i$, where $t_1, t_2, ..., t_n$ are nonnegative integers and $a_1, a_2, ..., a_n$ are given positive integers such that $\gcd(a_1, a_2, ..., a_n) = 1$. Finding the Frobenius number is called the Frobenius problem, the coin problem or the money changing problem(MCP). the Frobenius problem is not only interesting area for pure mathematicians but connected with graph theory in \cite{Heap (1965), Hujter (1987)}, the theory of computer science in \cite{Raczunas1996} as introduced in \cite{Owens2003}.
There is an explicit formula for the Frobenius number when there are only two relatively prime numbers in \cite{Sylvester1883}.
\\
\\
E.Curtis proved that the Frobenius number for three or more relatively prime numbers cannot be given by a finite set of polynomials in \cite{Curtis (1990)} and Ramirez-Alfonsin proved that the problem is NP-hard in \cite{Ramirez1996}. At the moment, there are only algorithmic methods to get the general formula for the Frobenius number of the set which has three or more relatively prime numbers in \cite{Beihoffer (2005), Bocker (2007)}. There are some recent results that the running time for the fastest algorithm is $O(a_1)$ with the residue table in memory in \cite{Brauer (1962)} and $O(na_1)$ with no additional memory requirements in \cite{Bocker (2007)}. Also, the research on limiting distribution in \cite{Shchur2009} and lower bound in \cite{Aliev (2007),Fel (2015)} of the Frobenius number was presented. In algebraic view, instead of finding the general formula for three or more relatively prime numbers, the formulae for special cases were found such as the Frobenius number of a set of integers in a geometric sequence in \cite{Ong2008}, a Pythagorean triples in \cite{Gil (2015)} and three consecutive squares or cubes in \cite{Lepilov (2015)}. Recently, the various ways to solve the Frobenius problem for numerical semigroups are suggested in \cite{Aliev (2009),Rosales2000,Rosales2009,Rosales2004}, etc. Especially, the method to compute the Apery set and to get the Frobenius number by using the Apery set is one of the efficient tool to solve the Frobenius problem of the numerical semigroups in \cite{Marquez (2015),Rosales2009,Ramirez2009}. For example, in the recent articles which presented the Frobenius problems for Fibonacci numerical semigroups in \cite{Marin (2007)}, Mersenne numerical semigroups and Thabit numerical semigroups in \cite{Rosales2015}, the method is used to get the Frobenius number.
\\
\\
The Frobenius problem in the numerical semigroups \\
$\big<\{3\cdot 2^{n+i} - 1| i \in \{0,1,\cdots \}\}\big>$ for $n \in \{0,1,\cdots \}$ was given in \cite{Rosales2015}. In \cite{Rosales2015}, the authors defined a Thabit number $3 \cdot 2^n - 1$ and Thabit numerical semigroups $T(n) = \big< \{3 \cdot 2^{n+i} - 1 | i \in \{0,1,\cdots\}\}\big>$ for a nonnegative integer $n$. And then, they defined minimal system of generators for $T(n)$ which is the smallest subset of $\big< \{3 \cdot 2^{n+i} - 1 | i \in \mathbb{N}\}\big>$ that equals $T(n)$. In \cite{Rosales2015}, the minimal system of generators for $T(n)$ is $\big< \{3 \cdot 2^{n+i} - 1 | i \in \{0,1,\cdots ,n+1\}\}\big>$ and it is unique. The embedding dimension is the cardinality of the minimal system of generators. By the minimality of the system $\big< \{3 \cdot 2^{n+i} - 1 | i \in \{0,1,\cdots ,n+1\}\}\big>$ for $T(n)$, the embedding dimension for $T(n)$ is $n+2$. For any set $S$ and $x \in S\texttt{\char`\\}\{0\}$, the Ap$\acute{e}$ry set was defined by $Ap(S,x) = \{s \in S | s - x \not\in S\}$. Let $s_i = 3 \cdot 2^{n+i} - 1$ for each nonnegative integer $i$. Then the Ap\'{e}ry set is defined by $Ap(T(n),s_0) = \{s \in T(n) | s-s_0 \not \in T(n) \}$ for $s_0$. In \cite{Rosales2015}, $Ap(T(n),s_0)$ was described explicitly which leads to solve the Frobenius problem.  Let $R(n)$ be the set of the sequences $(t_1,\cdots,t_{n+1}) \in \{0,1,2\}^{n+1}$ that satisfies the following conditions:
\begin{enumerate}
\item $t_{n+1} \in \{0,1\}$,
\item If $t_j = 2$, then $t_i = 0$ for all $i < j \leq n$,
\item If $t_n = 2$, then $t_{n+1} = 0$,
\item If $t_n = t_{n+1} = 1$, $t_i = 0$ for all $1\leq i <n$.
\end{enumerate}
Then \cite{Rosales2015} concludes that $Ap(T(n),s_0) = \{t_1 s_1 + \cdots + t_{n+1} s_{n+1} | (t_1,\cdots, t_{n+1}) \in R(n)\}$. The Frobenius number of the numerical semigroups was presented by $F(S) = \max(Ap(S,x)) - x$ in \cite{Rosales2009} and therefore the Frobenius number of Thabit numerical semigroups is $s_n + s_{n+1} - s_0 = 9\cdot 2^{2n} - 3\cdot 2^n - 1$.\\
\\
The main purpose of this paper is to extend the coefficients $3$ and $1$ in Thabit numerical semigroups in \cite{Rosales2015} to $2^k + 1$ and $2^k - 1$ for a positive integer $k$ and to define the minimal system of generators and the Ap\'{e}ry set. Also, we discuss the Frobenius problem in this numerical semigroups. More precisely, we solve the Frobenius problem for generalized Thabit numerical semigroups defined by $\{(2^k + 1)\cdot 2^{n+i} - (2^k - 1) | i \in \{0,1,\cdots \} \}$  for $n \in \{0,1,\cdots \}, k \in \{1,2,\cdots \}$. The paper in \cite{Rosales2015} is the case of $k = 1$ in this paper. 

This paper is organized as follows:
In Section 2, we introduce some facts that is required to understand this paper. 
In Section 3, we introduce the minimal system of generators and the embedding dimension for generalized Thabit numerical semigroups.
In Section 4, we introduce the progress to get the Ap\'{e}ry set and Frobenius number for generalized Thabit numerical semigroups.
In Section 5, we summarize the results for the Frobenius problem for generalized Thabit numerical semigroups.

\section{Preliminaries} \label{sec_Preliminaries}
Let $\mathbb{N}$ be the set of nonnegative integers.
\begin{definition} \cite{Rosales2015,Rosales2009}
A {\em numerical semigroup} is a subset $S$ of $\mathbb{N}$ that is closed under addition, contains $0$ and $\mathbb{N}\texttt{\char`\\}S$ is finite. Given a nonempty subset $A$ of $\mathbb{N}$ we will denote by $\big<A\big>$ the submonoid of $(\mathbb{N},+)$ generated by $A$, that is,\\
\begin{displaymath}
\big<A\big> = \{\lambda_1 a_1 + \cdots + \lambda_n a_n | n \in \mathbb{N}\texttt{\char`\\}\{0\}, a_i \in A, \lambda_i \in \mathbb{N} \textrm{ for all } i \in \{1,\cdots,n\}\}.
\end{displaymath}
\end{definition}
\begin{theorem}\cite{Rosales2015,Rosales2009}
$\big<A\big>$ is a numerical semigroup if and only if $\gcd(A) = 1$.
\end{theorem}
\begin{definition}
A positive integer $x$ is a {\em generalized Thabit number} if $x = (2^k + 1)\cdot 2^n - (2^k - 1)$ for some $n \in \mathbb{N}$ and $k \in \mathbb{N}\texttt{\char`\\}\{0\}$.
\end{definition}
\begin{definition}
A numerical semigroup $S$ is called a {\em generalized Thabit numerical semigroup} if there exists $n \in \mathbb{N}$ and $k \in \mathbb{N}\texttt{\char`\\}\{0\}$ such that $S = \big< \{(2^k + 1)\cdot 2^{n+i} - (2^k - 1) | i \in \mathbb{N}\}\big>$. We will denote by $GT(n,k)$ the generalized Thabit numerical semigroup $ \big< \{(2^k + 1)\cdot 2^{n+i} - (2^k - 1) | i \in \mathbb{N}\}\big>$.
\end{definition}
\begin{definition}\cite{Rosales2015,Rosales2009}
If $S$ is a numerical semigroup and $S = \big< A\big>$ then we say that $A$ is a {\em system of generators of} $S$. Moreover, if $S \neq \big< X\big>$ for all $X \not \subseteq A$, we say that $A$ is a {\em minimal system of generators of} $S$.
\end{definition}
\begin{theorem}\cite{Rosales2009}
Every numerical semigroup admits a unique minimal system of generators, which in addition is finite.
\end{theorem}
\begin{definition}\cite{Rosales2015,Rosales2009}
We call the cardinality of its minimal system of generators the {\em embedding dimension of} $S$, denoted by $e(S)$.
\end{definition}
\begin{definition}\cite{Rosales2015,Rosales2009}
We call the greatest integer that does not belong to $S$ the {\em Frobenius number of} $S$ denoted by $F(S)$ And the cardinality of $\mathbb{N}\texttt{\char`\\}S$ is called the {\em gender of} $S$ and denoted by $g(S)$.
\end{definition}
\section{The embedding dimension for $GT(n,k)$} \label{sec_The embedding dimension for $GT(n,k)$}
If $n \in \mathbb{N}$ and $k \in \mathbb{N}\texttt{\char`\\}\{0\}$, then $GT(n,k)$ is a submonoid of $(\mathbb{N},+)$. Moreover we have $\{(2^k + 1)\cdot 2^n -(2^k - 1), (2^k + 1)\cdot 2^{n+1} - (2^k - 1)\} \subseteq GT(n,k)$ and
\begin{align*}
\gcd((2^k + 1)\cdot 2^n & - (2^k - 1), (2^k + 1)\cdot 2^{n+1} - (2^k - 1)) \\
& = \gcd((2^k + 1)\cdot 2^n - (2^k - 1), 2^k - 1)\\
& = \gcd(2^{n+1}, 2^k - 1)\\
& = 1
\end{align*}
since $2^{n+1}$ has only even divisor and $2^k - 1$ is an odd number. Hence $\gcd(GT(n,k)) = 1$ and $GT(n,k)$ is a numerical semigroup.
\begin{lemma}\label{Lem_2a_2m_eq}
Let $A$ be a nonempty set of positive integers, $k \in \mathbb{N}\texttt{\char`\\}\{0\}$ and $M = \big<A\big>$. Then the following conditions are equivalent:
\begin{enumerate}
\item $2a+(2^k - 1) \in M$ for all $a \in A$,
\item $2m+(2^k - 1) \in M$ for all $m \in M\texttt{\char`\\}\{0\}$.
\end{enumerate}
\end{lemma}
{\it Proof.}\\
(1) $\Rightarrow$ (2). If $m \in M\texttt{\char`\\}\{0\}$, then there exist $a_1, \cdots, a_k \in A$ such that $m = a_1 + \cdots + a_k$($a_i = a_j$ for $i \neq j$ may be allowed). If $k=1$, then $m=a_1$ and therefore $2m+(2^k - 1) = 2a_1 + (2^k - 1) \in M$. If $k \geq 2$, then $2m+(2^k - 1) = 2(a_1 + \cdots + a_{k-1}) + 2a_k + (2^k - 1) \in M$, because $M$ is closed under addition.\\
(2) $\Rightarrow$ (1). Since all $a \in A$ can be represented by $a = a_1 = m \in M$ and hence $2a + (2^k - 1) = 2m + (2^k - 1) \in M$ for all $a \in A$.

\begin{proposition}
If $n \in \mathbb{N}$ and $k \in \mathbb{N}\texttt{\char`\\}\{0\}$, then $2t+(2^k - 1) \in GT(n,k)$ for all $t \in GT(n,k)\texttt{\char`\\}\{0\}$.
\end{proposition}
{\it Proof.}\\
Let $n \in \mathbb{N}$ and let $GT(n,k) = \big<\{(2^k + 1)\cdot 2^{n+i} - (2^k - 1)| i \in \mathbb{N}\}\big>$. Clearly $2\{(2^k + 1)\cdot 2^{n+i} - (2^k -1)\} + (2^k - 1) = (2^k + 1)\cdot 2^{n+i+1} - (2^k - 1) \in GT(n,k)$. From Lemma \ref{Lem_2a_2m_eq}, we obtain that $2t+(2^k - 1) \in GT(n,k)$ for all $t \in GT(n,k)\texttt{\char`\\}\{0\}$.

We need some preliminary results to find out the minimal system of generators of $GT(n,k)$.
\begin{lemma}\label{Lem_2s}
Let $n \in \mathbb{N}$, $k \in \mathbb{N}\texttt{\char`\\}\{0\}$, and
$S = \big<\{(2^k + 1)\cdot 2^{n+i} - (2^k - 1)| i \in \{0,1,\cdots,n+k\}\}\big>$. Then $2s+(2^k - 1) \in S$ for all $s \in S\texttt{\char`\\}\{0\}$.
\end{lemma}
{\it Proof.}\\
If $i \in \{0,1,\cdots,n+k-1\}$, then $2\{(2^k + 1)\cdot 2^{n+i} - (2^k - 1)\} + (2^k - 1) = (2^k + 1)\cdot 2^{n+i+1} - (2^k - 1) \in S$. Furthermore,
\begin{align*}
2\{(2^k & + 1)\cdot 2^{2n+k} - (2^k - 1)\} + (2^k - 1) \\
& = (2^k + 1)\cdot 2^{2n+k+1} - (2^k - 1) \\
& = \{(2^k + 1)\cdot 2^n - (2^k - 1)\}^2 + \{(2^k + 1)\cdot 2^{n+1} - (2^k - 1)\} \\
& \quad + (2^{k-1} - 1)\{(2^k + 1)2^{n+2} - (2^k - 1) \} + \{(2^k + 1)\cdot 2^{2n} - (2^k - 1)\} \\
& \quad + (2^{k-1} - 1)\{(2^k + 1)2^{2n+1} - (2^k - 1)\} \in S.
\end{align*} By using Lemma \ref{Lem_2a_2m_eq}, we obtain the desired result.

The next result gives a system of generators of $GT(n,k)$.
\begin{lemma}
If $n \in \mathbb{N}$ and $k \in \mathbb{N}\texttt{\char`\\}\{0\}$, then \[ GT(n,k) = \big<\{(2^k + 1)\cdot 2^{n+i} - (2^k - 1) | i \in \{0,1,\cdots, n+k\}\}\big>. \]
\end{lemma}
{\it Proof.}\\
Let $S = \big<\{(2^k + 1)\cdot 2^{n+i} - (2^k - 1) | i \in \{0,1,\cdots, n+k\}\}\big>$. Then $S \subseteq GT(n,k)$ is trivial. To show that
$S \supseteq GT(n,k)$ we need to check that $(2^k + 1)\cdot 2^{n+i} - (2^k - 1) \in S$ for all $i \in \mathbb{N}$. We use induction on $i$. For $i = 0$, the result is trivial.
Assume that the statement holds for $i$ and let us show it for $i+1$. Since $(2^k + 1)\cdot 2^{n+i+1} - (2^k - 1) = 2\{(2^k + 1)\cdot 2^{n+i} - (2^k - 1)\} + (2^k - 1)$ then, by Lemma \ref{Lem_2s}, we get that $(2^k + 1)\cdot 2^{n+i+1} - (2^k - 1) \in S$. Hence by the principle of mathematical induction the result holds for all $i$ in nonnegative integers.

The next result shows that $(2^k + 1)\cdot 2^{2n+k} - (2^k - 1)$ belongs to the minimal system of generators of $GT(n,k)$ if $k \leq n$.\\
\begin{lemma}\label{Lem_minimal_k<n}
If $n \in \mathbb{N}$ and $k \in \mathbb{N}\texttt{\char`\\}\{0\}$, then $(2^k + 1)\cdot 2^{2n+k} - (2^k - 1) \not \in \big<\{(2^k + 1)\cdot 2^{n+i} - (2^k - 1) | i \in \{0,1,\cdots, n+k-1\}\}\big>$ for $k \leq n$.
\end{lemma}
{\it Proof.}\\
Let us suppose on the conversely that $(2^k + 1)\cdot 2^{2n+k} - (2^k - 1) \in \big<\{(2^k + 1)\cdot 2^{n+i} - (2^k - 1) | i \in \{0,1,\cdots, n+k-1\}\}\big>$. Then there exist $a_0, \cdots, a_{n+k-1} \in \mathbb{N}$ such that
\begin{align*}
& (2^k + 1)\cdot 2^{2n+k} - (2^k - 1) \\
& = a_0 \{(2^k + 1)\cdot 2^n - (2^k - 1)\} + \cdots + a_{n+k-1} \{(2^k + 1)\cdot 2^{2n+k-1} - (2^k - 1)\}\\
& = (2^k + 1)(a_0 2^n + \cdots + a_{n+k-1} 2^{2n+k-1}) - (2^k - 1)(a_0 + \cdots + a_{n+k-1})
\end{align*}
and consequently $a_0 + \cdots + a_{n+k-1} \equiv 1 \mod{(2^k + 1)\cdot 2^{n}}$. Hence $a_0 + \cdots + a_{n+k-1} = 1 + t \cdot (2^k + 1) \cdot 2^n$ for some $t \in \mathbb{N}$. Besides, it is clear that $t \neq 0$(If $t = 0$, $(2^k + 1)\cdot 2^{2n+k} - (2^k - 1) = a_0 \{(2^k + 1)\cdot 2^n - (2^k - 1)\} + \cdots + a_{n+k-1} \{(2^k + 1)\cdot 2^{2n+k-1} - (2^k - 1)\} \leq (2^k + 1)\cdot 2^{2n+k-1} - (2^k - 1)$ and it is a contradiction.) and so $a_0 + \cdots + a_{n+k-1} \geq 1 + (2^k + 1)\cdot 2^n$.
\begin{align*}
a_0 \{(2^k & + 1)\cdot 2^n - (2^k - 1)\} + \cdots + a_{n+k-1} \{(2^k + 1)\cdot 2^{2n+k-1} - (2^k - 1)\}\\
& \geq (a_0 + \cdots + a_{n+k-1})\{(2^k + 1) \cdot 2^n - (2^k - 1)\}\\
& \geq \{1 + (2^k + 1)\cdot 2^n\}\{(2^k + 1) \cdot 2^n - (2^k - 1)\}\\
& = (2^k + 1)^2 2^{2n} - (2^k - 2)(2^k + 1)2^n - (2^k - 1)\\
& =  2^{2n+2k} + 2^{2n+k} + 2^{2n+k} + 2^{2n} - 2^{n+2k} + 2^{n+k} + 2^n + 2^n - (2^k - 1)\\
& = \{(2^k + 1)\cdot 2^{2n+k} - (2^k - 1) \} + 2^{2n+k} + 2^{2n} - 2^{n+2k} + 2^{n+k} + 2^n + 2^n\\
& > \{(2^k + 1)\cdot 2^{2n+k} - (2^k - 1) \} \textrm{ if } k \leq n \textrm{ since } 2^{2n+k} \geq 2^{n+2k}.
\end{align*}
Therefore, it completes the proof.

And we suggest the conclusion for a minimal system of generators of $GT(n,k)$ in the following Theorem \ref{Thm_minimal}.
\begin{theorem}\label{Thm_minimal}
Let $n \in \mathbb{N}$, $k \in \mathbb{N}\texttt{\char`\\}\{0\}$ and
\begin{displaymath}
\delta = \left\{ \begin{array}{ll}
1 & \text{if} ~~ n = 0,\\
k & \text{if} ~~ n \neq 0, k \leq n,\\
k - 1 & \text{if} ~~ n \neq 0, k > n.
\end{array} \right.
\end{displaymath}
Then $\{(2^k + 1)\cdot 2^{n+i} - (2^k - 1)| i \in \{0,\cdots, n + \delta\}\}$ is a minimal system of generators.
\end{theorem}
{\it Proof.}\\
Let us suppose that there exist $a_0, \cdots, a_{n+k-1} \in \mathbb{N}$ such that
\begin{align*}
& (2^k + 1)\cdot 2^{2n+k} - (2^k - 1) \\
& = a_0 \{(2^k + 1)\cdot 2^n - (2^k - 1)\} + \cdots + a_{n+k-1} \{(2^k + 1)\cdot 2^{2n+k-1} - (2^k - 1)\}\\
& = (2^k + 1)(a_0 2^n + \cdots + a_{n+k-1} 2^{2n+k-1}) - (2^k - 1)(a_0 + \cdots + a_{n+k-1})
\end{align*}
as Lemma \ref{Lem_minimal_k<n} and let us suppose that $a_0 + \cdots + a_{n+k-1} = 1 + 2 \cdot (2^k + 1) \cdot 2^n$. In other words, $t = 2$ in this case. Note that the following calculation is similar as the progress to get Lemma \ref{Lem_minimal_k<n}:
\begin{align*}
& a_0 \{(2^k + 1)\cdot 2^n - (2^k - 1)\} + \cdots + a_{n+k-1} \{(2^k + 1)\cdot 2^{2n+k-1} - (2^k - 1)\}\\
& \geq (a_0 + \cdots + a_{n+k-1})\{(2^k + 1) \cdot 2^n - (2^k - 1)\}\\
& \geq \{1 + (2^k + 1)\cdot 2^{n+1}\}\{(2^k + 1) \cdot 2^n - (2^k - 1)\}\\
& = (2^k + 1)^2 2^{2n+1} - (2^k - 1)(2^k + 1)\cdot 2^{n+1} + (2^k + 1) \cdot 2^n - (2^k - 1)\\
& = 2^{2n+2k+1} + 2^{2n+k+1} + 2^{2n+k+1} + 2^{2n+1} - 2^{n+2k+1} + 2^{n+1} + 2^{n+k} + 2^n \\
& \quad - (2^k - 1)\\
& = 2^{2n+2k} + 2^{2n+k} + 2^{2n+2k} + 2^{2n+k} + 2^{2n+k+1} + 2^{2n+1} - 2^{n+2k+1} + 2^{n+1}\\
& \quad + 2^{n+k} + 2^n - (2^k - 1)\\
& = \{(2^k + 1)\cdot 2^{2n+k} - (2^k - 1)\} + 2^{2n+2k} + 2^{2n+k} + 2^{2n+k+1} + 2^{2n+1}\\
& \quad - 2^{n+2k+1} + 2^{n+1} + 2^{n+k} + 2^n\\
& > \{(2^k + 1)\cdot 2^{2n+k} - (2^k - 1)\} \textrm{ if } n \neq 0 \textrm{ since } 2^{2n+2k} \geq 2^{n+2k+1}.
\end{align*}
By summarizing the above results we can notice that $(2^k + 1)\cdot 2^{2n+k} - (2^k - 1) \not \in \big<\{(2^k + 1)\cdot 2^{n+i} - (2^k - 1) | i \in \{0,1,\cdots, n+k-1\}\}\big>$ if ($n \neq 0$ and $t \geq 2$). Therefore, we only have to check the remaining cases : (1) $t = 1$ or (2) $n = 0$. But we can easily notice that $(2^k + 1)\cdot 2^0 - (2^k - 1) = 2$ for all $k$ and it implies that for the case of $n = 0$, a minimal system of generators of $GT(n,k)$ is $\{(2^k + 1)\cdot 2^{i} - (2^k - 1) | i \in \{0,1\}\}$ and hence we will only check for the case of $t = 1$.\\
Let $t = 1$. In other words, $a_0 + \cdots + a_{n+k-1} = 1 + (2^k + 1)\cdot 2^{n}$. Then, we can set the system of equations such that
\begin{gather}
a_0 2^n + \cdots + a_{n+k-1}2^{2n+k-1} = 2^{2n+k} + (2^k - 1)\cdot 2^n, \label{eq_1}\\
a_0 + \cdots + a_{n+k-1} = 1 + (2^k + 1)\cdot 2^n. \label{eq_2}
\end{gather}
Divide the equation \ref{eq_1} by $2^n$, we obtain 
\begin{gather}
a_0 + a_1 \cdot 2 + \cdots + a_{n+k-1} \cdot 2^{n+k-1} = 2^{n+k} + (2^k - 1).\label{eq_3}
\end{gather}
Note that the equation (\ref{eq_3}) $-$ (\ref{eq_2}) is : \\
$a_1 + \cdots + a_{n+k-1} \cdot (2^{n+k-1} - 1) = 2^{n+k} + (2^k - 1) - (2^k + 1) \cdot 2^n - 1 = 2^k - 2^n - 2$.\\
We can assume that $k > n$ since the minimality of\\
$\big<\{(2^k + 1)\cdot 2^{n+i} - (2^k - 1) | i \in \{0,1,\cdots, n+k\}\}\big>$ is guaranteed for $k \leq n$ by Lemma \ref{Lem_minimal_k<n}. Note that $a_k = a_{k+1} = \cdots = a_{n+k-1} = 0$ since $2^k - 2^n - 2 < 2^k - 1$ and hence $a_1 + a_2 (2^2 - 1) + \cdots + a_{k-1}(2^{k-1} - 1) = 2^k - 2^n - 2$. We can find a solution for any $n$ and $k$. It is $a_1 = 2^k - 2^n - 2$ and $a_2 = \cdots = a_{k-1} = 0$. By adjusting $a_0$ to satisfy (2), we can find a solution to represent $(2^k + 1)\cdot 2^{2n+k} - (2^k - 1)$ by $\big<\{(2^k + 1)\cdot 2^{n+i} - (2^k - 1) | i \in \{0,1,\cdots, n+k-1\}\}\big>$. It is $(a_0 , a_1 , a_2 , \cdots , a_{n+k-1}) = (2^{n+k} + 2^{n+1} - 2^k + 3, 2^k - 2^n - 2, 0,\cdots , 0)$.\\
Note that the solution for $(a_0 , a_1 , a_2 , \cdots , a_{n+k-1})$ is not unique.\\
Since $2^k - 2^n - 2 \geq 0$ for any $k > n \geq 1$, the remaining step is to check the minimality of $\big<\{(2^k + 1)\cdot 2^{n+i} - (2^k - 1) | i \in \{0,1,\cdots, n+k-1\}\}\big>$ for $k > n$.\\
Let us Assume that there exists $\alpha \leq k - 2$ such that $\big<\{(2^k + 1)\cdot 2^{n+i} - (2^k - 1) | i \in \{0,1,\cdots, n+ \alpha \}\}\big> = GT(n,k)$. In other words, the following system of equations has a solution.
\begin{gather}
a_0 2^n + \cdots + a_{n+\alpha }2^{2n+\alpha } = 2^{2n+\alpha + 1} + (2^k - 1)\cdot 2^n, \label{eq_4}\\
a_0 + \cdots + a_{n + \alpha } = 1 + (2^k + 1)\cdot 2^n. \label{eq_5}
\end{gather}
Similar to the progress to getting the equation \ref{eq_3}, we can get the equation by dividing the equation (\ref{eq_4}) by $2^n$,
\begin{gather}
a_0 + a_1 \cdot 2 + \cdots + a_{n + \alpha } \cdot 2^{n + \alpha } = 2^{n + \alpha + 1} + (2^k - 1). \label{eq_6}
\end{gather}
And the equation (\ref{eq_6}) $-$ (\ref{eq_5}) is :
\begin{align*}
a_1 + \cdots + a_{n+\alpha } \cdot (2^{n + \alpha } - 1)
& = 2^{n+\alpha + 1} + (2^k - 1) - (2^k + 1) \cdot 2^n - 1 \\
& = 2^{n+\alpha + 1} - 2^{n+k} + 2^k - 2^n - 2 \\
& = (2^{n+\alpha + 1} - 2^{n+k}) + 2^k - 2^n - 2 \\
& \leq -2^{n+k-1} + 2^k - 2^n - 2 \\
& \leq -2^n - 2 \\
& < 0  \qquad(\text{since } \alpha \leq k-2).
\end{align*}
This implies that to guarantee the minimality, the number of elements in generator is at least $n + k$ for $k > n$. It completes the proof.

By Theorem \ref{Thm_minimal}, we can identify the embedding dimension of $GT(n,k)$ for all $n \in \mathbb{N}$, and $k \in \mathbb{N}\texttt{\char`\\}\{0\}$.
\begin{corollary}
Let $n \in \mathbb{N}$, $k \in \mathbb{N}\texttt{\char`\\}\{0\}$ and let $GT(n,k)$ be a generalized Thabit numerical semigroup associated to $n$ and $k$. Then $e(GT(n,k)) = n + \delta + 1$.
\end{corollary}
Gathering all this information we obtain that for each integer $t$ greater than or equal to $2$ there exists a generalized Thabit numerical semigroup $GT(n,k)$ with embedding dimension $t$.
\begin{example}
$GT(0,3) = \big<\{(2^3 + 1) \cdot 2^0 - (2^3 - 1), (2^3 + 1) \cdot 2^1 - (2^3 - 1) \}\big> = \big<\{2, 11\}\big>$ is a generalized Thabit numerical semigroup with embedding dimension $0 + 1 + 1 = 2$.
\end{example}

\begin{example}
$GT(3,2)= \big< \{(2^2 + 1) \cdot 2^3 - (2^2 - 1), (2^2 + 1) \cdot 2^4 - (2^2 - 1), (2^2 + 1) \cdot 2^5 - (2^2 - 1),(2^2 + 1) \cdot 2^6 - (2^2 - 1), (2^2 + 1) \cdot 2^7 - (2^2 - 1), (2^2 + 1) \cdot 2^8 - (2^2 - 1) \}\big> = \big<\{37,77,157,317,637,1277\}\big>$ is a generalized Thabit numerical semigroup with embedding dimension $3 + 2 + 1 = 6$.
\end{example}

\begin{example}
$GT(1,3) = \big< \{(2^3 + 1) \cdot 2^1 - (2^3 - 1), (2^3 + 1) \cdot 2^2 - (2^3 - 1), (2^3 + 1) \cdot 2^3 - (2^3 - 1), (2^3 + 1) \cdot 2^4 - (2^3 - 1) \}\big> = \big< \{11,29,65,137\}\big>$ is a generalized Thabit numerical semigroup with embedding dimension $1+(3-1)+1 = 4$.
\end{example}
\newpage
\section{The Ap\'{e}ry set for $GT(n,k)$}\label{sec_The Apery set for $GT(n,k)$}
Let $S$ be a numerical semigroup and let $x \in S \texttt{\char`\\}\{0\}$. The Ap\'{e}ry set of $x$ in $S$ is defined as $Ap(S,x) = \{s \in S | s-x \not \in S \}$ in \cite{Rosales2009}. From \cite{Rosales2009}, we easily deduce the following Lemma.
\begin{lemma}
Let $S$ be a numerical semigroup and let $x \in S \texttt{\char`\\}\{0\}$. Then $Ap(S,x)$ has cardinality equal to $x$. Moreover $Ap(S,x) = \{w(0),w(1),\cdots,\\
w(x-1)\}$ where $w(i)$ is the least element of $S$ congruent with $i$ modulo $x$ for all $i \in \{0,\cdots, x-1\}$.
\end{lemma}
{\it Proof.}\\
It is in the proof of Lemma 2.4 in \cite{Rosales2009}.

\begin{example}
Let $S = \left<\{7,11,13\}\right>$. Then
\[ S = \{ 0,7,11,13,14,18,20,21,22, 24, 25, 26, 27,28,29,31,\rightarrow\} \]
where the symbol $\rightarrow$ means that every integer greater than $31$ belongs to the set. Hence $Ap(S,7) = \{0,11,13,22,24,26,37\}$.
\end{example}
The next result is due to Selmer in \cite{Selmer1977} (also, see \cite{Rosales2009}) and can be used to solve the Frobenius problem for a generalized Thabit numerical semigroup from its Ap\'{e}ry set. We provide the proof in more detail.
\begin{lemma}\label{lem_F_g}
Let $S$ be a numerical semigroup and let $x \in S\texttt{\char`\\}\{0\}$. Then,
\begin{enumerate}
\item $F(S) = \max(Ap(S,x)) - x$\\
\item $g(S) = \frac{1}{x} (\sum_{w \in Ap(S,x)} w) - \frac{x-1}{2}$
\end{enumerate}
\end{lemma}
{\it Proof.}\\ (1) By the definition of the Ap\'{e}ry set, $\max(Ap(S,x)) - x \not \in S$. Hence $F(S) \geq \max(Ap(S,x)) - x$.
Let us assume that $\max(Ap(S,x)) - x + t \not \in S$ for some $t \in \mathbb{N}\texttt{\char`\\}\{0\}$. In other words, $F(S) > \max(Ap(S,x)) - x$. Note that $\max(Ap(S,x)) + t \in S$ is trivial since $\max(Ap(S,x)) + t = w(i) + t_1 x$ for some $t_1 \in \mathbb{N}\texttt{\char`\\}\{0\}$. Since $\max(Ap(S,x)) - x + t = (\max(Ap(S,x)) + t) - x$, $\max(Ap(S,x)) + t \in Ap(S,x)$ by assumption. But it is a contradiction. Hence, $F(S) = \max(Ap(S,x)) - x$.\\
(2) It is in the proof of Proposition 2.12 in \cite{Rosales2009}.

From now on we will denote by $s_i$ the elements $(2^k + 1)\cdot 2^{n+i} - (2^k - 1)$ for each $i \in \{0,1,\cdots,n+\delta \}$. Thus with this notation we have that $\{s_0,s_1,\cdots,s_{n+\delta}\}$ is the minimal system of generators of $GT(n,k)$.
\begin{lemma}\label{Lem_push_coefficients_1}
Let $n \in \mathbb{N}$ and $k \in \mathbb{N}\texttt{\char`\\}\{0\}$. Then:
\begin{enumerate}
\item If $0<i\leq j<n+\delta $ then $s_i + 2s_j = 2s_{i-1} + s_{j+1}$.\\
\item If $0<i\leq n+\delta $ then
\end{enumerate}
\begin{displaymath}
s_i + 2s_{n+\delta } =\left\{ \begin{array}{ll}
2s_{i-1} + s_1 + (2^k + 1)s_0 & \text{if} ~~ n = 0,\\
\\
2s_{i-1} + {s_0}^2 + s_1 + (2^{k-1} - 1)s_2 \\
+ s_{n + k - 1} + (2^{k-1} -1)s_{n+k } & \text{if} ~~ n \neq 0, k \leq n,\\
\\
2s_{i-1} + (2^{n+k} + 2^{n+1} - 2^k + 3)s_0 \\
+ (2^k - 2^n - 2)s_1 & \text{if} ~~ n \neq 0, k > n.
\end{array} \right.
\end{displaymath}
\end{lemma}
{\it Proof.}\\ (1) If $0<i\leq j<n+\delta $, then we have that
\begin{align*}
s_i + 2s_j & = (2^k + 1)\cdot 2^{n+i} - (2^k - 1) + 2\{(2^k + 1)\cdot 2^{n+j} - (2^k - 1)\} \\
& = 2\{(2^k + 1)\cdot 2^{n+i-1} - (2^k - 1)\} + (2^k + 1)\cdot 2^{n+j+1} -(2^k - 1) \\
& = 2s_{i-1} + s_{j+1}.
\end{align*}
(2) For convenience, if $0 < i \leq n + \delta$, we can consider three cases :
\begin{enumerate}
\item If $n = 0$, then
\begin{align*}
s_i + & 2s_{n + \delta} \\
&= s_i + 2s_1 \\
&= \{(2^k + 1)\cdot 2^i - (2^k - 1)\} + 2\{(2^k + 1)\cdot 2^1 - (2^k - 1)\} \\
&= 2\{(2^k + 1) \cdot 2^{i-1} - (2^k - 1)\} + \{(2^k + 1) \cdot 2^2 - (2^k - 1)\} \\
&= 2s_{i-1} \\
& \quad + \{(2^k + 1)\cdot 2^1 - (2^k - 1)\} + (2^k + 1)\{(2^k + 1)\cdot 2^0 - (2^k - 1)\} \\
& = 2s_{i-1} + s_1 + (2^k + 1)s_0.
\end{align*}
\item If $n \neq 0, k \leq n$, then
\begin{align*}
s_i + & 2s_{n+\delta } \\
& = (2^k + 1)\cdot 2^{n+i} - (2^k - 1) + 2\{(2^k + 1)\cdot 2^{2n+k} - (2^k - 1)\} \\
& = (2^k + 1)\cdot 2^{n+i} - 2(2^k - 1) + (2^k + 1)\cdot 2^{2n+k+1} - (2^k - 1) \\
& = 2\{(2^k + 1)\cdot 2^{n+i-1} - (2^k - 1)\} + \{(2^k + 1) \cdot 2^n - (2^k - 1)\}^2 \\
& \quad + \{(2^k + 1) \cdot 2^{n+1} - (2^k - 1)\}\\
& \quad + (2^{k-1} - 1)\{(2^k + 1)\cdot 2^{n+2} - (2^k - 1)\}\\
& \quad + \{(2^k + 1)\cdot 2^{2n} - (2^k - 1)\}\\
& \quad + (2^{k-1} -1)\{(2^k + 1)\cdot 2^{2n+1} - (2^k - 1)\} \\
& = 2s_{i-1} + {s_0}^2 + s_1 + (2^{k-1} - 1)s_2 + s_n + (2^{k-1} -1)s_{n+1}.
\end{align*}
\item If $n \neq 0, k > n$, then
\begin{align*}
s_i + & 2s_{n+\delta} \\
& = (2^k + 1)\cdot 2^{n+i} - (2^k - 1) + 2\{(2^k + 1)\cdot 2^{2n+k-1} - (2^k - 1)\} \\
& = (2^k + 1)\cdot 2^{n+i} - 2(2^k - 1) + (2^k + 1)\cdot 2^{2n+k} - {2^k - 1} \\
& = 2\{(2^k + 1)\cdot 2^{n+i-1} - (2^k - 1)\} \\
& \quad + (2^{n+k} + 2^{n+1} - 2^k + 3)\{(2^k + 1)\cdot 2^n - (2^k - 1)\} \\
& \quad + (2^k - 2^n - 2)\{(2^k + 1)\cdot 2^{n+1} - (2^k - 1)\} \\
& = 2s_{i-1} + (2^{n+k} + 2^{n+1} - 2^k + 3)s_0 + (2^k - 2^n - 2)s_1.
\end{align*}
\end{enumerate}

By Lemma \ref{Lem_push_coefficients_1}, we can consider the set of coefficients $(t_1,\cdots, t_{n+\delta})$ such that the expressions $t_1 s_1 + \cdots + t_{n + \delta} s_{n + \delta}$ represent all elements in \\
$Ap(GT(n,k),s_0)$. We will conduct step by step approach to establish the set of coefficients $(t_1,\cdots,t_{n+\delta})$. First, we get the set of coefficients $(t_1,\cdots , t_{n+\delta})$ such that $t_1 s_1 + \cdots + t_{n + \delta} s_{n + \delta}$ contains all elements in $Ap(GT(n,k),s_0)$ but may not be equal.  We can obtain the set by the following Lemma.
\begin{lemma} \label{Lem_A}
Let $A(n,\delta)$ be the set of $(t_1, \cdots, t_{n+\delta}) \in \{0,1,2\}^{n+\delta}$ such that $t_{n+\delta} = 0$ or $1$ and if $t_j = 2$, then $t_i = 0$ for all $i < j$. Then $Ap(GT(n,k),s_0) \subseteq \{t_1 s_1 + \cdots + t_{n + \delta} s_{n + \delta} | (t_1,\cdots,t_{n+\delta}) \in A(n,\delta)\}$ unless $n = 1$ and $k = 2$.
\end{lemma}
{\it Proof.}\\
It is same as the proof of Lemma 10 in \cite{Rosales2015} except the case of ($n = 1$, $k = 2$) and I will show that $2s_{n+\delta} \not \in Ap(GT(n,k),s_0)$ unless ($n = 1$, $k = 2$) in Corollary \ref{cor_n1_k2}.

\begin{lemma}\label{Lem_x-1}
Let $n \in \mathbb{N}$ and $k \in \mathbb{N}\texttt{\char`\\}\{0\}$. If $x \in GT(n,k)$ and $x \not\equiv0 \mod{s_0}$, $x - (2^k - 1) \in GT(n,k)$.
\end{lemma}
{\it Proof.}\\
If $x \in GT(n,k)$ then there exist $a_0,\cdots,a_{n+\delta} \in \mathbb{N}$ such that $x = a_0 s_0 + \cdots + a_{n+\delta} s_{n+\delta}$. On the other hand, if $x \not\equiv 0 \mod{s_0}$ then there exists $i \in \{1,\cdots, n+\delta\}$ such that $a_i \neq 0$. Hence
\begin{align*}
& x - (2^k - 1)\\
& =  a_0 s_0 + \cdots + (a_i - 1)s_i + \cdots + a_{n+\delta}s_{n+\delta} + (2^k + 1)\cdot 2^{n+i} - 2(2^k - 1)\\
& = a_0 s_0 + \cdots + (a_i - 1)s_i + \cdots + a_{n+\delta}s_{n+\delta} + 2\{(2^k + 1)\cdot 2^{n+i-1} - (2^k - 1)\}\\
& = a_0 s_0 + \cdots + (a_{i-1} + 2) s_{i-1} + (a_i - 1)s_i + \cdots + a_{n+\delta} s_{n+\delta} \\
& \quad \in GT(n,k).
\end{align*}

Note that on contrary to this Lemma, $x \in T(n)$ implies that $x - 1 \in T(n)$ in \cite{Rosales2015}.
Lemma \ref{Lem_x-1} is very important since $\gcd(2^k - 1, s_0) = \gcd(2^k - 1, (2^k + 1)\cdot 2^n - (2^k - 1)) = \gcd(2^k - 1, (2^k + 1)\cdot 2^n) = \gcd(2^k - 1, 2^{n+1}) = 1$, for any $x > (2^k - 1)\{(2^k + 1)\cdot 2^n - 2^k\}$, the set $\{x - i(2^k - 1)| i \in \{0,1,\cdots,(2^k + 1)\cdot 2^n - 2^k\}\}$ is a complete system of residues modulo $s_0$ and hence we can get the following Lemma.
\begin{lemma}\label{Lem_maxAp_1}
If $n \in \mathbb{N}$ and $k \in \mathbb{N}\texttt{\char`\\}\{0\}$ then
\begin{gather*}
w(s_0 - (2^k - 1)) = \max(Ap(GT(n,k),s_0)).
\end{gather*}
\end{lemma}
From this Lemma and the fact that $s_i \equiv (2^i - 1)(2^k - 1)\mod{s_0}$, we can prove the following Corollary previously announced in Lemma \ref{Lem_A}.
\begin{corollary}\label{cor_n1_k2}
If $n \in \mathbb{N}$ and $k \in \mathbb{N}\texttt{\char`\\}\{0\}$ then $2s_{n+\delta} \not \in Ap(GT(n,k),s_0)$ unless $n=1$ and $k = 2$.
\end{corollary}
Note that
\begin{displaymath}
\delta = \left\{ \begin{array}{ll}
1 & \text{if} ~~ n = 0,\\
k & \text{if} ~~ n \neq 0, k \leq n,\\
k - 1 & \text{if} ~~ n \neq 0, k > n.
\end{array} \right.
\end{displaymath}
{\it Proof.}\\
We already know in the proof of Lemma \ref{Lem_A} that $s_i + 2s_{n+\delta} \not \in Ap(GT(n,k),s_0)$ for all $i \in \mathbb{N}$. Therefore, if $2s_{n+\delta} \in Ap(GT(n,k),s_0)$, it should be a maximum element. In other words, $2s_{n+\delta} = max(Ap(GT(n,k),s_0))$. Let us consider the following three cases :
\begin{enumerate}
\item If $n = 0$, then
\begin{align*}
2s_{n+\delta} - s_0 & = 2s_1 - s_0 \\
& = 2\{(2^k + 1)\cdot 2 - (2^k - 1)\} - 2 \\
& = 2\cdot 2^k + 6 - 2 \\
& = 2\cdot 2^k + 4 \\
& = 2\cdot (2^{k-1} + 2) \\
& = (2^{k-1} + 2)s_0 \in GT(0,k).
\end{align*}
\item If $n \neq 0$ and $k \leq n$, then since $s_i \equiv (2^i - 1)(2^k - 1)\mod{s_0}$ and hence $2s_{n+k} \equiv (2^{n+k+1} - 2)(2^k - 1)\mod{s_0}$. Note that
\begin{align*}
(2^{n+k} + 2^n - 2^k)(2^k - 1) & \equiv \{2(2^{n+k} + 2^n - 2^k) + 1\}(2^k - 1)\\
& \equiv -(2^k - 1) \mod{s_0}
\end{align*}
and
\begin{align*}
(2^{n+k} + 2^n - 2^k)  < 2^{n+k+1} - 2 
 < \{2(2^{n+k} + 2^n - 2^k) + 1\}.
\end{align*} This implies that $(2^{n+k+1} - 2)(2^k - 1) \not \equiv -(2^k - 1) \mod{s_0}$ and hence $2s_{n+\delta} \not \in Ap(GT(n,k),s_0)$ if $n \neq 0$, $k \leq n$.\\
\item If $n \neq 0$ and $k > n$, then
$2s_{n+k-1} \equiv (2^{n+k} - 2)(2^k - 1)\mod{s_0}$ and
\begin{align*}
(2^{n+k} + 2^n - 2^k) & < 2^{n+k} - 2\\
& < \{2(2^{n+k} + 2^n - 2^k) + 1\}
\end{align*}

is satisfied unless $n = 1$ and $k = 2$. Hence $2s_{n+\delta} \not \in Ap(GT(n,k),s_0)$ if $n \neq 0$, $k > n$ and $(n,k) \neq (1,2)$. If $(n,k) = (1,2)$, $2^{n+k} + 2^n - 2^k = 2^{n+k} - 2$ is satisfied and
\begin{align*}
2s_{n+k-1} - s_0
& = 2s_2 - s_0 \\
& = 2\{(2^2 + 1)\cdot 2^3 - (2^2 - 1)\} - \{(2^2 + 1)\cdot 2 - (2^2 - 1)\} \\
& = 67.
\end{align*}
Note that it cannot be represented by $\big<\{s_0,s_1,s_2\}\big> = \big<\{7,17,37\}\big>$.
\end{enumerate}

We will check the following four cases to find the Frobenius number and the Ap\'{e}ry set for $GT(n,k)$:
\begin{enumerate}
\item $n = 0$,\\
\item $2 \leq k \leq n$,\\
\item $n \neq 0, 2 \neq k > n$,\\
\item $n = 1$ and $k = 2$.
\end{enumerate}
Note that the case of $k = 1$ is previously solved in \cite{Rosales2015}.
\subsection{The case of $n = 0$}
We will state the Frobenius number and the Ap\'{e}ry set for the case in the following Lemma.
\begin{lemma}\label{Lem_maxAp_n=0}
If $n = 0$ and $k \in \mathbb{N}\texttt{\char`\\}\{0\}$, $\max(Ap(GT(0,k), s_0)) = s_1$ and hence $F(GT(0,k)) = s_1 - s_0 = s_1 - 2$. Moreover, $Ap(GT(0,k),s_0) = \{0, s_1\}$.
\end{lemma}
{\it Proof.}\\
It is trivial since $2 \in GT(0,k)$ and $(2^k + 1) \cdot 2^i - (2^k - 1)$ is odd number for all $i \geq 1$.

\subsection{The case of $2 \leq k \leq n$}
Let us consider $s_i + s_{n+k} \in GT(n,k)$ for $i \neq 0$.
\begin{align*}
s_i + s_{n+k} - s_0
& = \{(2^k + 1) \cdot 2^{n+i} - (2^k - 1)\} + \{(2^k + 1) \cdot 2^{2n+k} - (2^k - 1)\} \\
& \quad - \{(2^k + 1) \cdot 2^{n} - (2^k - 1)\}\\
& = (2^k + 1) \cdot 2^{n}(2^i + 2^{n+k} - 1) - (2^k - 1).
\end{align*}
Then we can set the system of equations:
\begin{gather*}
\alpha_0 + \cdots + 2^{n+k} \alpha_{n+k} = 2^i + 2^{n+k} - 1 + (2^k - 1) \beta, \\
\alpha_0 + \cdots + \alpha_{n+k} = 1 + (2^k + 1)\cdot 2^n \beta.
\end{gather*}
Therefore $\alpha_1 + (2^2 - 1)\alpha_2 + \cdots + (2^{n+k} - 1)\alpha_{n+k} = 2^i + 2^{n+k} - 2 + \beta (2^k - 1 - 2^{n+k} - 2^n)$.
If $\beta = 0$, $s_i + s_{n+k} - s_0 \leq s_{n+k}$ and it is a contradiction. Hence, if there is a solution $(\alpha_0, \cdots, \alpha_{n+k})$ when $\beta = 1$, this implies $s_i + s_{n+k} \not \in Ap(GT(n,k), s_0)$.
Let $\beta = 1$. Then
\begin{align*}
2^i + 2^{n+k} - 2 + & \beta (2^k - 1 - 2^{n+k} - 2^n) \\
& = 2^i + 2^{n+k} - 2 + 2^k - 1 - 2^{n+k} - 2^n \\
& = 2^i + 2^k - 2^n - 3.
\end{align*}
Hence $s_i + s_{n+k} \not \in Ap(GT(n,k),s_0)$ if (i) $k = n$ and $i \geq 2$ or (ii) $i \geq n$.\\
To improve the result of (i), let us consider
\begin{align*}
2s_1 + s_{n+k} - s_0
& = 2s_1 + s_{2n} - s_0\\
& = (2^n + 1)\cdot 2^{n+2} - 2(2^n - 1) + (2^n + 1)\cdot 2^{3n}\\
& \quad - (2^n - 1) - (2^n + 1)\cdot 2^n + (2^n - 1)\\
& = (2^n + 1)\cdot 2^n (2^2 + 2^{2n} - 1) - 2(2^n - 1).
\end{align*}
We can establish the system of equations
\begin{gather}
\alpha_0 + 2\alpha_1 + \cdots + 2^{2n} \alpha_{2n} = 2^2 + 2^{2n} - 1 + (2^n - 1)\beta \label{eq_7}\\
\alpha_0 + \cdots + \alpha_{2n} = 2 + (2^n + 1)2^n \beta \label{eq_8}
\end{gather}
Hence, the equation (\ref{eq_7}) $-$ (\ref{eq_8}) is
\begin{align*}
\alpha_1 + (2^2 - 1)\alpha_2 + \cdots + (2^{2n} - 1)\alpha_{2n} & = 1 + 2^{2n} + \beta (2^n - 1 - 2^{2n} - 2^n) \\
& = 1 + 2^{2n} - \beta (1 + 2^{2n}).
\end{align*}
Note that there is a solution $\alpha_1 = \cdots = \alpha_{2n} = 0$ and $\alpha_0 = 2^{2n} + 2^n + 2$ for $\beta = 1$. Hence $2s_1 + s_{n+k} \not \in Ap(GT(n,k),s_0)$.
Then, let us consider $s_1 + s_{n+k} - s_0$ and we can set the system of equations
\begin{gather*}
\alpha_0 + 2\alpha_1 + \cdots + 2^{2n} \alpha_{2n} = 2 + 2^{2n} - 1 + (2^n - 1) \beta \\
\alpha_0 + \cdots + \alpha_{2n} = 1 + (2^n + 1)2^n \beta.
\end{gather*}
Then
\begin{align*}
\alpha_1 + (2^2 - 1)\alpha_2 + \cdots + (2^{2n} - 1)\alpha_{2n} & = 2^{2n} + \beta (2^n - 1 - 2^{2n} - 2^n)\\
& = 2^{2n} - \beta (1 + 2^{2n}).
\end{align*}
And this implies that there is no solution if $\beta \geq 1$.
Hence, we can get the following Lemma.
\begin{lemma}\label{Lem_maxAp_k=n}
If $k = n \geq 1$, $\max(Ap(GT(n,k),s_0)) = s_1 + s_{n+k}$ and $F(GT(n,n)) = s_1 + s_{n+n} - s_0 = (2^n + 1)\cdot 2^n (2^{2n} + 1) - (2^n - 1)$.
\end{lemma}
Actually, for the one who want to prove only Lemma \ref{Lem_maxAp_k=n}, the Lemma can be proved more easily by using Lemma \ref{Lem_maxAp_1} since there is only one element whose residue is $-(2^n - 1)$ modulo $s_0$ in $\{s_2 + s_{n+n}, 2s_1 + s_{n+n}, s_1 + s_{n+n}\}$.\\
To improve the result of (ii) in previous page, let us consider $2s_{n-1} + s_{n+k} - s_0$ and we can establish the system of equations
\begin{gather*}
\alpha_0 + 2\alpha_1 + \cdots + 2^{n+k} \alpha_{n+k} = 2^n + 2^{n+k} - 1 + (2^k - 1)\beta , \\
\alpha_0 + \cdots + \alpha_{n+k} = 2 + (2^k + 1)\cdot 2^n \beta .
\end{gather*}
Then, we can obtain
\begin{align*}
\alpha_1 + (2^2 - 1)\alpha_2 & + \cdots + (2^{n+k} - 1) \alpha_{2n+k} \\
& = 2^n + 2^{n+k} - 3 + (2^k - 1 - 2^{n+k} - 2^n)\beta
\end{align*}
and therefore $2s_{n-1} + s_{n+k} \not \in Ap(GT(n,k),s_0)$ if $2 \leq k < n$.\\
Let us consider
$s_{n-1} + s_{n+k} - s_0$ and we can set the system of equations
\begin{gather*}
\alpha_0 + \cdots + 2^{n+k} \alpha_{n+k} = 2^{n-1} + 2^{n+k} - 1 + (2^k - 1) \beta, \\
\alpha_0 + \cdots + \alpha_{n+k} = 1 + (2^k + 1)2^n \beta .
\end{gather*}
Then
\begin{align*}
\alpha_1 + (2^2 - 1)\alpha_2 & + \cdots + (2^{n+k} - 1) \alpha_{2n+k} \\
& = 2^{n-1} + 2^{n+k} - 2 + (2^k - 1 - 2^{n+k} - 2^n) \beta.
\end{align*}
If $\beta = 1$,
\begin{align*}
2^{n-1} + 2^{n+k} - 2 + & (2^k - 1 - 2^{n+k} - 2^n) \beta \\
& = 2^{n-1} + 2^{n+k} - 2 + (2^k - 1 - 2^{n+k} - 2^n) \\
& = -2^{n-1} + 2^k - 3
\end{align*}
and by $k < n$, there is no solution. Hence we can establish the following Lemma.
\begin{lemma}
If $2 \leq k < n$, $ s_{n-1} + s_{n+k} \leq \max(Ap(GT(n,k),s_0)) < 2s_{n-1} + s_{n+k}$.
\end{lemma}
Let us consider $t_1 s_1 + t_2 s_2 + \cdots + t_{n-2} s_{n-2} + s_{n-1} + s_{n+k} - s_0$ where\\
$(t_1,\cdots,t_{n-2},1, 0, \cdots, 0, 1) \in A(n,\delta )$. Then
\begin{align*}
t_1 s_1 & + t_2 s_2 + \cdots + t_{n-2} s_{n-2} + s_{n-1} + s_{n+k} - s_0\\
& = t_1\{(2^k + 1)\cdot 2^{n+1} - (2^k - 1)\} + t_2\{(2^k + 1)\cdot 2^{n+2} - (2^k - 1)\} + \cdots\\
& \quad + t_{n-2}\{(2^k + 1)\cdot 2^{2n-2} - (2^k - 1)\} + \{(2^k + 1)\cdot 2^{2n-1} - (2^k - 1)\} \\
& \quad + \{(2^k + 1)\cdot 2^{2n+k} - (2^k - 1)\} - (2^k + 1)\cdot 2^{n} + (2^k - 1)\\
& = (2^k + 1)\cdot 2^n (2t_1 + 2^2 t_2 + \cdots + 2^{n-2} t_{n-2} + 2^{n-1} + 2^{n+k} - 1) \\
& \quad - (t_1 + \cdots + t_{n-2} + 1)(2^k - 1).
\end{align*}
And we can obtain the system of equations
\begin{align*}
\alpha_0 + & 2\alpha_1 + \cdots + 2^{n+k} \alpha_{n+k} \\
& = (2t_1 + 2^2 t_2 + \cdots + 2^{n-2} t_{n-2} + 2^{n-1} + 2^{n+k} - 1) + (2^k - 1)\beta ,\\
\alpha_0 + & \cdots + \alpha_{n+k} = (t_1 + \cdots + t_{n-2} + 1) + (2^k + 1)\cdot 2^n \beta .
\end{align*}
Then if $\beta = 1$,
\begin{align*}
\alpha_1 + & (2^2 - 1)\alpha_2 + \cdots + (2^{n+k} - 1) \alpha_{n+k}\\
& = t_1 + (2^2 - 1)t_2 + \cdots + (2^{n-2} - 1)t_{n-2} + 2^{n-1} + 2^{n+k} - 2 \\
& \quad + \beta (2^k - 1 - 2^{n+k} - 2^n)\\
& = t_1 + (2^2 - 1)t_2 + \cdots + (2^{n-2} - 1)t_{n-2} - 2^{n-1} + 2^k - 3.
\end{align*}
Hence, we can establish the following Corollary.
\begin{corollary}
If $k = n-1$, then $\max(Ap(GT(n,k),s_0)) = 2s_1 + s_{n-1} + s_{n+k}$ and $F(GT(n,n-1)) = 2s_1 + s_{n-1} + s_{2n-1} - s_0 = (2^{n-1} + 1)\cdot 2^n (3 + 2^{n-1} + 2^{2n-1}) - 3(2^{n-1} - 1)$. Also, if $k = n - 2$, then $\max(Ap(GT(n,k),s_0)) = s_2 + s_{n-2} + s_{n-1} + s_{n+k}$.
\end{corollary}
Note that if
\begin{align*}
t_1 + (2^2 - 1)t_2 +& \cdots + (2^{n-2} - 1)t_{n-2} - 2^{n-1} + 2^k - 3\\
& = t_1 + (2^2 - 1)t_2 + \cdots + (2^{n-2} - 1)t_{n-2} - (2^{n-1} - 2^k + 3)\\
& \geq 0,
\end{align*}
we can conclude that $t_1s_1 + \cdots + t_{n-2}s_{n-2} + s_{n-1} + s_{n+k} \not\in Ap(GT(n,k),s_0)$ and finally, we can obtain $\max(Ap(GT(n,k),s_0))$ for all $2 \leq k < n$.
\begin{lemma} \label{Lem_maxAp_2<k<n}
Let $2\leq k < n$ and let the coefficients $(t_1,\cdots,t_{n-2}, 1, 0,\cdots, 0, 1) \in A(n,\delta)$ such that
\begin{align*}
t_1 + (2^2 - 1)t_2 + \cdots + (2^{n-2} - 1)t_{n-2} = 2^{n-1} - 2^k + 2
\end{align*} which is presented in above. Then
\begin{align*}
\max(Ap(GT(n,k),s_0)) = t_1 s_1 + \cdots + t_{n-2} s_{n-2} + s_{n-1} + s_{n+k}
\end{align*}
and
\begin{align*}
F(GT(n,k)) = {} & \max(Ap(GT(n,k),s_0)) - s_0 \\
= {} & t_1 s_1 + \cdots + t_{n-2} s_{n-2} + s_{n-1} + s_{n+k} - s_0.
\end{align*}
\end{lemma}
\begin{example}
$GT(5,3) \\= \big<\{281,569,1145,2297,4601,9209,18425,36857,73721\}\big>$ by Theorem \ref{Thm_minimal}. Since $2^{n-1} - 2^k + 2 = (2^2 - 1) \cdot 1 + (2^3 - 1) \cdot 1 = 10$, we can conclude that $\max(Ap(GT(5,3),281) = s_2 + s_3 + s_4 + s_8 = 81764$ and $F(GT(5,3)) = \max(Ap(GT(5,3),281)) - 281 = 81764 - 281 = 81483$.
\end{example}
\begin{corollary}
$F(\left<\{5\cdot 2^{n+i} - 3 | i \in \mathbb{N}, n \geq 2\}\right>) = 85\cdot 2^{2n-2} - 5\cdot 2^n - 9$.
\end{corollary}
{\it Proof.}\\ 
It can be obtained by Lemma \ref{Lem_maxAp_2<k<n} with fixed $k = 2$. You can easily find the solution $(t_1,\cdots,t_{n-2}, 1, 0, \cdots, 0, 1) \in A(n,\delta)$ of the equation such that
\begin{gather*}
t_1 + (2^2 - 1)t_2 + \cdots + (2^{n-2} - 1)t_{n-2} = 2^{n-1} - 2^2 + 2
\end{gather*}
and the solution is $t_{n-2} = 2$ and $t_i = 0$ for all $i \leq n-3$. Hence
\begin{align*}
& F(\left<\{5\cdot 2^{n+i} - 3 | i \in \mathbb{N}, n \geq 2\}\right>) \\
& = 2s_{n-2} + s_{n-1} + s_{n+k} - s_0 \\
& = 2\cdot (5\cdot 2^{2n-2} - 3) + 5\cdot 2^{2n-1} - 3 + 5\cdot 2^{2n+2} - 3 - (5\cdot 2^{n} - 3) \\
& = 5\cdot 2^{2n-2} + 5\cdot 2^{2n+2} - 5\cdot 2^n - 9 \\
& = 85\cdot 2^{2n-2} - 5\cdot 2^n - 9.
\end{align*}
We introduce an interesting Lemma which is used to more efficient computation to get the Frobenius number in special cases.
\begin{lemma}\label{Lem_Inequality_MaxAp_1}
If $2 \leq k < n$, then $s_k + \cdots + s_{n-1} + s_{n+k} \in Ap(GT(n,k),s_0)$.
\end{lemma}
{\it Proof.}\\
Let $k = n - 1 - \gamma$ such that $\gamma \geq 1$. Then $2^{n-1} - 2^k + 3 = 2^{n-1-\gamma}(2^{\gamma} - 1) + 3$ and by $2^{\gamma} - 1 = 1 + \cdots + 2^{\gamma - 1}$, we can establish the condition of $(t_1,\cdots, t_{n-2})$ such that $t_1 s_1 + \cdots + t_{n-2} s_{n-2} + s_{n-1} + s_{n+k} \not \in Ap(GT(n,k),s_0)$ for $2 \leq k < n$ by the following inequality :
\begin{gather}\label{eq_9}
t_1 + (2^2 - 1)t_2 + \cdots + (2^{n-2} - 1)t_{n-2} \geq 2^{n-1-\gamma}(1 + \cdots + 2^{\gamma - 1}) + 3.
\end{gather}
By setting $t_{n-1-\gamma} = \cdots = t_{n-2} = 1$ and $t_i = 0$ for $i \neq n-1-\gamma, \cdots, n-2$, the LHS in the equation \ref{eq_9} is smaller than RHS in the equation \ref{eq_9} by $\gamma + 3$. Therefore we can notice that for $(t_1,\cdots,t_{n-2},1, 0, \cdots, 0, 1) \in A(n,\delta ), t_1s_1 + \cdots + t_{n-2}s_{n-2} + s_{n-1} + s_{n+k} \in Ap(GT(n,k),s_0)$ if $t_k = \cdots = t_{n-2} = 1$ is satisfied.

By using Lemma \ref{Lem_maxAp_2<k<n} and Lemma \ref{Lem_Inequality_MaxAp_1}, we can obtain the following Corollary and for sufficiently small $n$, we can use the following Corollary to calculate $\max(Ap(GT(n,k),s_0))$ more easily.

\begin{corollary}
Let $2\leq k < n$ and $\gamma + 3 = n + 2 - k \leq 2^k -1$. In other words, $n \leq 2^k + k - 3$. And let $(t_1,\cdots,t_{k-1},1,\cdots,1(=t_{n-1}),0,\cdots,0,1) \in A(n,\delta)$ be an element of $A(n,\delta)$ which satisfies the following inequality :
\begin{gather}
t_1 + (2^2 - 1)t_2 + \cdots + (2^{k-1} - 1)t_{k-1} = n + 1 - k.
\end{gather}
Then
\begin{gather*}
\max(Ap(GT(n,k),s_0)) = t_1 s_1 + \cdots + t_{k-1} s_{k-1} + s_k + \cdots + s_{n-1} + s_{n+k}.
\end{gather*}
\end{corollary}
\begin{example}
Since $2 \leq 3 < 5$ and $5 \leq 2^3 + 3 - 3$, we can get\\
$\max(Ap(GT(5,3),281))$ more easily. Since $(2^2 - 1) \cdot 1 = 5 + 1 - 3 = 3$, we can conclude that \[\max(Ap (GT(5,3) ,281)) = s_2 + s_3 + s_4 + s_8 = 81764. \]
\end{example}
\begin{example}
We can also consider $GT(7,3) =\\ \big<\{1145,2297,4601,9209,18425,36857,73721,147449,294905,589817,1179641\}\big>$. Since $7 \leq 2^3 + 3 - 3$, we can easily notice that $2 + 1 \cdot (2^2 - 1) = 5 = n + 1 - k$ leads to $\max(Ap(GT(7,3), 1145)) = 2s_1 + s_2 + s_3 + s_4 + s_5 + s_6 + s_{10} = 1327048$ and $F(GT(7,3)) = \max(Ap(GT(7,3), 1145)) - 1145 = 1325903$.
\end{example}
Let $R_1 (n,\delta)$ be the set of the sequences $(t_1,\cdots,t_{n+\delta}) \in A(n,\delta)$ that satisfies the following conditions if $t_{n+k} = 1$:
\begin{enumerate}
\item $t_n = \cdots = t_{n+k-1} = 0$,\\
\item If $t_{n-1} = 1$, $t_1 + (2^2 - 1)t_2 + \cdots + (2^{n-2} - 1)t_{n-2} \leq 2^{n-1} - 2^k + 2$.
\end{enumerate}
Then, we obtain the following Lemma:
\begin{lemma}\label{Lem_Ap_k<n}
Let $2 \leq k < n$. Then
\begin{gather*}
Ap(GT(n,k), s_0) = \{t_1 s_1 + \cdots + t_{n + \delta} s_{n + \delta} | (t_1,\cdots,t_{n+\delta}) \in R_1 (n,\delta)\}.
\end{gather*}
\end{lemma}
{\it Proof.}\\
Note that $Ap(GT(n,k), s_0) \subseteq \{t_1 s_1 + \cdots + t_{n + \delta} s_{n + \delta} | (t_1,\cdots,t_{n+\delta}) \in R_1 (n,\delta)\}$ and $\#\{t_1 s_1 + \cdots + t_{n + \delta} s_{n + \delta} | (t_1,\cdots,t_{n+\delta}) \in R_1(n,\delta)\} \leq \# R_1(n,\delta)$. Then if $\# R_1(n,\delta)$ has the cardinality $s_0 = (2^k + 1) \cdot 2^n - (2^k - 1)$, it suffices the proof. We classify the cases to obtain the cardinality as the following :
\begin{enumerate}
\item If $t_{n+\delta} = 0$, then it can be classified again to two cases as the following: \label{2<k<n_case_1}
\begin{enumerate}
\item Let $2 \not \in \{t_1, \cdots, t_{n+\delta - 1}\}$. Then
$\#\{(t_1,\cdots,t_{n+\delta}) \in R_1(n,\delta)\} = 2^{n+\delta-1}$ since $t_i \in \{0,1\}$ for all $1 \leq i \leq n+\delta-1$.
\item Let $2 \in \{t_1, \cdots, t_{n+\delta - 1}\}$ \label{2<k<n_case_1_2}. If $t_i = 2$ for some $i \in \{1,\cdots,n+\delta - 1\}$, then $t_j = 0$ for all $j < i$ and $t_j \in \{0,1\}$ for all $i < j \leq n+\delta - 1$. Hence, $\#R_1(n,\delta)$ in this case is $2^{n+\delta - i - 1}$. So, we use the summation to find $\#R_1(n,\delta)$ for (\ref{2<k<n_case_1_2}) : $\sum_{i=1}^{n+\delta - 1} 2^{n+\delta - i - 1} = 2^{n+\delta - 1} - 1$.
\end{enumerate}
Hence, $\#R_1(n,\delta) = 2^{n+\delta} - 1$ in case \ref{2<k<n_case_1}.\\
\item If $t_{n+\delta} = 1$, then it can be classified again to two cases as the following:
\begin{enumerate}
\item Let $t_{n-1} = 1$. Then
we can find the unique solution $(t_1,\cdots,t_{n-2})$ for each $0 \leq x \leq 2^{n-1} - 2$ : $(t_1,\cdots,t_{n-2}) \in A(n,\delta)$ such that $t_1 + (2^2 - 1)t_2 + \cdots + (2^{n-2} - 1)t_{n-2} = x$. Hence, the number of the solution $(t_1, \cdots, t_{n-2})$ satisfies the inequality $(t_1,\cdots,t_{n-2}) \in A(n,\delta)$ such that $t_1 + (2^2 - 1)t_2 + \cdots + (2^{n-2} - 1)t_{n-2} \leq 2^{n-1} - 2^{\delta} + 2$ is $2^{n-1} - 2^{\delta} + 3$.
\item Let $t_{n-1} = 0$.  \label{2<k<n_case_2_2} then it can be classified again to two cases as the following:
\begin{enumerate}
\item If $2 \not \in \{t_1,\cdots, t_{n-2}\}$ then  $\#R_1(n,\delta) = 2^{n-2}$.
\item If $2 \in \{t_1,\cdots,t_{n-2}\}$ then $\#R_1(n,\delta) = \sum_{i=1}^{n - 2} 2^{n-2-i} = 2^{n-2} - 1$.
\end{enumerate}
Note that we can verify $\#R_1(n,\delta) = 2^{n-1} - 1$ in case (\ref{2<k<n_case_2_2}).
\end{enumerate}
\end{enumerate}
Hence, we can conclude that $\#R_1(n,\delta) = 2^{n+\delta} - 1 + 2^{n-1} - 2^{\delta} + 3 + 2^{n-1} - 1 = 2^{n+\delta} + 2^n - 2^{\delta} + 1 = (2^{\delta} + 1)\cdot 2^n - (2^{\delta} - 1) = (2^{k} + 1)\cdot 2^n - (2^k - 1) = s_0$ since $\delta = k$ for $k \leq n$.

\begin{example}
Let $(n,k)=(3,2)$. Since the solutions of the inequality $t_1 + (2^2 - 1)t_2 + \cdots + (2^{n-2} - 1)t_{n-2} \leq 2^{3-1} - 2^2 + 2$ are $(0,0,\cdots,0),(1,0,\cdots,0)$ and $(2,0,\cdots,0)$, we can obtain
\begin{align*}
Ap(&GT(3,2),37) \\
& = \{0,s_1,2s_1,s_2,s_1+s_2,\\
& \quad 2s_1+s_2,2s_2,s_3,s_1+s_3,2s_1+s_3,\\
& \quad s_2+s_3,s_1+s_2+s_3,2s_1+s_2+s_3,2s_2+s_3,2s_3,\\
& \quad s_4,s_1+s_4,2s_1+s_4,s_2+s_4,s_1+s_2+s_4,\\
& \quad 2s_1+s_2+s_4,2s_2+s_4,s_3+s_4,s_1+s_3+s_4,2s_1+s_3+s_4,\\
& \quad s_2+s_3+s_4,s_1+s_2+s_3+s_4,2s_1+s_2+s_3+s_4,2s_2+s_3+s_4,2s_3+s_4,\\
& \quad 2s_4,s_5,s_1+s_5,2s_1+s_5,s_2+s_5,\\
& \quad s_1+s_2+s_5, 2s_1+s_2+s_5\}.
\end{align*}
Note that $\# Ap(GT(3,2),37) = 37$ and $\max(Ap(GT(3,2),37)) = 2s_1 + s_2 + s_5$.
\end{example}
We already know from Lemma \ref{Lem_maxAp_k=n} that $\max(Ap(GT(n,k),s_0)) = s_1 + s_{n+k}$ if $k = n$ and $F(GT(n,n)) = s_1 + s_{n+n} - s_0 = (2^n + 1)\cdot 2^n (2^{2n} + 1) - (2^n - 1)$. Then $R_2 (n,\delta)$ can be defined by the set of the sequences $(t_1,\cdots,t_{n+\delta}) \in A(n,\delta)$ which satisfies the condition that if $t_{2n} = 1$, $t_1 = 1$ and $t_2 = \cdots = t_{2n-1} = 0$. Then we can obtain the following Lemma:
\begin{lemma}\label{Lem_Ap_k=n}
Let $k = n \geq 1$. Then
\begin{align*}
Ap(GT(n,k), s_0) = \{t_1 s_1 + \cdots + t_{n + \delta} s_{n + \delta} | (t_1,\cdots,t_{n+\delta}) \in R_2 (n,\delta)\}.
\end{align*}
\end{lemma}
{\it Proof.}\\
We can classify the cases by the following:
\begin{enumerate}
\item If $t_{2n} = 1$, then we can choose $t_1 \in \{0,1\}$ so there are 2 cases. \label{k=n_case_1} \\
\item If $t_{2n} = 0$, then it can be classified again to two cases as the following:
\begin{enumerate} \label{k=n_case_2}
\item Let $2 \not \in \{t_1,\cdots, t_{2n-1}\}$ Then we can choose the coefficients \\
$(t_1,\cdots,t_{2n-1})$ and each coefficient has the condition $t_i \in \{0,1\}$ for all $i$. Hence there are $2^{2n-1}$ cases.
\item Let $2 \in \{t_1,\cdots, t_{2n-1}\}$. Then if $t_i = 2$, we can choose the coefficients $(t_{i+1},\cdots,t_{2n-1})$ and each coefficient has the condition $t_j \in \{0,1\}$ for all $i < j \leq 2n-1$. Hence there are $2^{2n-1-i}$ cases for each $i$. Then we can express all the cases by the summation : $\sum_{i=1}^{2n-1} 2^{2n-1-i} = \sum_{i=0}^{2n-2} 2^i = 2^{2n-1} - 1$.
\end{enumerate}
\end{enumerate}
Therefore, we can obtain $2 + 2^{2n-1} + 2^{2n-1} - 1 = 2^{2n} + 1$ by case \ref{k=n_case_1} and case \ref{k=n_case_2}. Note that $2^{2n} + 1 = (2^n + 1)\cdot 2^n - (2^n - 1) = s_0$ and it completes the proof.

\begin{example}
Let $(n,k)=(2,2)$. Then we can obtain
\begin{align*}
Ap(&GT(2,2),17) \\
& = \{0,s_1,2s_1,s_2,s_1+s_2,\\
& \quad 2s_1+s_2,2s_2,s_3,s_1+s_3,2s_1+s_3,\\
& \quad s_2+s_3,s_1+s_2+s_3,2s_1+s_2+s_3,2s_2+s_3,2s_3,\\
& \quad s_4,s_1+s_4\}.
\end{align*}
Note that $\# Ap(GT(2,2),17) = 17$ and $\max(Ap(GT(2,2),17)) = s_1 + s_4$.
\end{example}
\subsection{The case of $n \neq 0, 2 \neq k > n$}
Let us consider
\begin{align*}
t_1 s_1 + \cdots + & t_{n+k-1} s_{n+k-1} - s_0 \\
& = t_1\{(2^k + 1) \cdot 2^{n+1} - (2^k - 1)\} + \cdots \\
& \quad + t_{n+k-1}\{(2^k + 1) \cdot 2^{2n+k-1} - (2^k - 1)\} - s_0\\
& = (2^k + 1)(2t_1 + 2^2 t_2 + \cdots + 2^{n+k-1} t_{n+k-1} - 1) \\
& \quad - (2^k - 1)(t_1 + \cdots + t_{n+k-1} - 1).
\end{align*}
Then we can set the system of equations\\
\begin{align*}
\alpha_0 + 2\alpha_1 + & \cdots + 2^{n+k-1} \alpha_{n+k-1} \\
& = (2t_1 + 2^2 t_2 + \cdots + 2^{n+k-1} t_{n+k-1} - 1) + (2^k - 1)\beta,\\
\alpha_0 + \alpha_1 + & \cdots + \alpha_{n+k-1} \\
& = t_1 + \cdots + t_{n+k-1} - 1 + (2^k + 1)\cdot 2^n \beta.
\end{align*}
It implies that
\begin{align*}
\alpha_1 + & (2^2 - 1)\alpha_2 + \cdots + (2^{n+k-1} - 1)\alpha_{n+k-1} \\
& = t_1 + (2^2 - 1)t_2 + \cdots + (2^{n+k-1} - 1)t_{n+k-1} + (2^k - 1 - 2^{n+k} - 2^n)\beta.
\end{align*}
And then if
\begin{align*}
\alpha_1 + & (2^2 - 1)\alpha_2 + \cdots + (2^{n+k-1} - 1)\alpha_{n+k-1} \\
& = t_1 + (2^2 - 1)t_2 + \cdots + (2^{n+k-1} - 1)t_{n+k-1} + (2^k - 1 - 2^{n+k} - 2^n) \\
& \geq 0,
\end{align*}
then $t_1 s_1 + \cdots + t_{n+k-1} s_{n+k-1} \not \in Ap(GT(n,k),s_0)$. Note that $2^{n+k} + 2^n - 2^k + 1 = (2^{n+k-1} - 1) + (2^{n+k-1} + 2^n - 2^k + 2) < 2\cdot (2^{n+k-1} - 1)$ if $n \neq 0, 2 \neq k > n$. Hence, we can find the $\max(Ap(GT(n,k),s_0))$ for all $n \neq 0, 2 \neq k > n$.
\begin{lemma}
Let $n \neq 0, 2 \neq k > n$ and let the coefficients\\
 $(t_1,\cdots,t_{n+k-2},1)\in A(n,\delta)$ such that
\begin{align*}
t_1 + (2^2 - 1)t_2 + \cdots + (2^{n+k-2} - 1)t_{n+k-2} =  2^{n+k-1} + 2^n - 2^k + 1.
\end{align*}
Then
\begin{gather*}
\max(Ap(GT(n,k),s_0)) = t_1 s_1 + \cdots + t_{n+k-2} s_{n+k-2} + s_{n+k-1}
\end{gather*}
and
\begin{align*}
F(GT(n,k)) & = \max(Ap(GT(n,k),s_0)) - s_0 \\
& =  t_1 s_1 + \cdots + t_{n+k-2} s_{n+k-2} + s_{n+k-1} - s_0.
\end{align*}
\end{lemma}

 Note that
\begin{align*}
2^{n+k-1} + & 2^n - 2^k + 1 \\
& = (2^{n-1} - 1) 2^k + (2^n + 1) \\
& = (1 + 2 + \cdots + 2^{n-2})2^k + (2^n + 1) \\
& = (2^k - 1) + \cdots + (2^{n+k-2} - 1) + (n - 1) + (2^n + 1) \\
& = (2^k - 1) + \cdots + (2^{n+k-2} - 1) + (2^n + n).
\end{align*}
Also, we can easily deduce the fact that $2^n + n < 2^{n+1} - 1 \leq 2^k - 1$ for $n \neq 1$ and $2^n + n = 2^{n+1} - 1 < 2^k - 1$ for $n = 1$. Therefore, we can reduce an effort to get $\max(Ap(GT(n,k),s_0))$ for $n \neq 0, 2 \neq k > n$.
\begin{lemma}\label{Lem_maxAp_k>n}
Let $n \neq 0, 2 \neq k > n$ and let the coefficients \\
$(t_1,\cdots,t_{k-1},1,\cdots,1) \in A(n,\delta)$ such that
\begin{align*}
t_1 + (2^2 - 1)t_2 + \cdots + (2^{k-1} - 1)t_{k-1} = 2^n + n.
\end{align*}
Then
\begin{align*}
\max(& Ap(GT(n,k), s_0)) \\
& = t_1 s_1 + \cdots + t_{k-1} s_{k-1} + s_{k} + s_{k+1} + \cdots + s_{n+k-2} + s_{n+k-1}
\end{align*} and
\begin{align*}
F(GT & (n,k)) \\
& = \max(Ap(GT(n,k),s_0)) - s_0 \\
& = t_1 s_1 + \cdots + t_{k-1} s_{k-1} + s_{k} + s_{k+1} + \cdots + s_{n+k-2} + s_{n+k-1} - s_0.
\end{align*}
\end{lemma}
\begin{example}
We have that $GT(2,3) = \big<\{ 29,65,137,281,569\} \big>$ and $2^n + n = 2^2 + 2 = 6$. So, $t_2 = 2$ and $t_i = 0$ for all $i < k$ and $i \neq 2$. Then you can get $\max(Ap(GT(2,3),s_0)) = 2s_2 + s_3 + s_4 = 1124$ and $F(GT(2,3)) = \max(Ap(GT(2,3),s_0)) - s_0 = 1124 - 29 = 1095$ by Lemma \ref{Lem_maxAp_k>n}.
\end{example}
$R_3(n,\delta)$ can be defined by the set of the sequences $(t_1,\cdots,t_{n+\delta}) \in A(n,\delta)$ which satisfies the condition that if $t_k = t_{k+1} = \cdots = t_{n+k} = 1$,
\begin{gather*}
t_1 + (2^2 - 1)t_2 + \cdots + (2^{k-1} - 1)t_{k-1} \leq 2^n + n.
\end{gather*}
Then we can obtain the following Lemma.
\begin{lemma}\label{Lem_Ap_k>n}
Let $n \neq 0, 2 \neq k > n$. Then
\begin{gather*}
Ap(GT(n,k),s_0) = \{t_1 s_1 + \cdots + t_{n + \delta} s_{n + \delta} | (t_1,\cdots,t_{n+\delta}) \in R_3(n,\delta)\}.
\end{gather*}
\end{lemma}
{\it Proof.}\\
Note that $Ap(GT(n,k), s_0) \subseteq \{t_1 s_1 + \cdots + t_{n + \delta} s_{n + \delta} | (t_1,\cdots,t_{n+\delta}) \in R_3 (n,\delta)\}$ and $\#\{t_1 s_1 + \cdots + t_{n + \delta} s_{n + \delta} | (t_1,\cdots,t_{n+\delta}) \in R_3(n,\delta)\} \leq \# R_3(n,\delta)$. Then if $\# R_3(n,\delta)$ has the cardinality $s_0 = (2^k + 1) \cdot 2^n - (2^k - 1)$, it suffices the proof. We can classify the cases to verify the cardinality as the following:
\begin{enumerate}
\item Let $(t_k,\cdots, t_{n+k-1}) \in \{0,1\}^n$ and there exists at least one \\
$t_i \in \{t_k,\cdots,t_{n+k-1}\}$ such that $t_i = 0$. Then there are $2^n - 1$ cases to choose $\{t_k,\cdots,t_{n+k-1}\}$ and to choose $\{t_1,\cdots,t_{k-1}\}$, it can be classified again to two cases as the following:
\begin{enumerate}
\item Let $2 \not \in \{t_1,\cdots,t_{k-1}\}$. Then there are $2^{k-1}$ cases to choose $\{t_1,\cdots,t_{k-1}\}$.
\item Let $2 \in \{t_1,\cdots, t_{k-1}\}$. Then if $t_i = 2$, there are $2^{k-i-1}$ cases to choose $\{t_{i+1},\cdots,t_{k-1}\}$. Hence, $\sum_{i=1}^{k-1} 2^{k-i-1} = 2^{k-1} - 1$.
\end{enumerate}
Therefore, we can check that case 1 have $(2^{k-1} + 2^{k-1} - 1)(2^n - 1) = (2^k - 1)(2^n - 1) = 2^{n+k} - 2^n - 2^k + 1$ cases.
\item Let $t_k = t_{k+1} = \cdots = t_{n+k-1} = 1$. Then we can obtain the unique solution $(t_1,\cdots,t_{k-1})$ for each $0 \leq x \leq 2^{n} + n$ such that $t_1 + (2^2 - 1)t_2 + \cdots + (2^{k-1} - 1)t_{k-1} = x$. Hence, the number of the solutions $(t_1, \cdots, t_{k-1})$ satisfy the inequality $t_1 + (2^2 - 1)t_2 + \cdots + (2^{k-1} - 1)t_{k-1} \leq 2^{n} + n$ is $2^n + n + 1$.
\item Let $2 \in \{t_k,\cdots, t_{n+k-2}\}$. If $t_i = 2$, there are $2^{n+k-i-1} - 1$ cases to choose $\{t_{i+1},\cdots,t_{n+k-1}\}$ since at least one $t_j \in \{t_{i+1},\cdots,t_{n+k-1}\}$ have to be $0$. Then we can obtain the summation : $\sum_{i=k}^{n+k-2} (2^{n+k-i-1} - 1) = \sum_{i=1}^{n-1} 2^i - (n - 1) = 2^n - 2 - (n - 1) = 2^n - n - 1$.
\end{enumerate}
You can verify easily that the equation
\begin{align*}
(2^{n+k} - 2^n - & 2^k + 1) + (2^n + n + 1) + (2^n - n - 1) \\
& = 2^{n+k} + 2^n - 2^k + 1 \\
& = 2^{n+k} + 2^n - (2^k - 1) \\
& = s_0
\end{align*}
 holds and it completes the proof.

\begin{example}
Let $(n,k)=(2,3)$. Since the solutions of the inequality $t_1 + (2^2 - 1)t_2 + \cdots + (2^{k-1} - 1)t_{k-1} \leq 2^2 + 2$ are $(0,0,\cdots,0),(1,0,\cdots,0),(2,0,\cdots,0),\\
(0,1,0,\cdots,0),(1,1,0,\cdots,0),(2,1,0,\cdots,0)$ and $(0,2,0,\cdots,0)$, we can obtain
\begin{align*}
Ap(&GT(2,3),29) \\
& = \{0,s_1,2s_1,s_2,s_1+s_2,\\
& \quad 2s_1+s_2,2s_2,s_3,s_1+s_3,2s_1+s_3,\\
& \quad s_2+s_3,s_1+s_2+s_3,2s_1+s_2+s_3,2s_2+s_3,2s_3,\\
& \quad s_4,s_1+s_4,2s_1+s_4,s_2+s_4,s_1+s_2+s_4,\\
& \quad 2s_1+s_2+s_4,2s_2+s_4,s_3+s_4,s_1+s_3+s_4,2s_1+s_3+s_4,\\
& \quad s_2+s_3+s_4,s_1+s_2+s_3+s_4,2s_1+s_2+s_3+s_4,2s_2+s_3+s_4\}.
\end{align*}
Note that $\# Ap(GT(2,3),29) = 29$ and $\max(Ap(GT(2,3),29)) = 2s_2+s_3+s_4$.
\end{example}

\subsection{The case of $n = 1$ and $k = 2$}
The Frobenius number and the Ap\'{e}ry set can be verified easily as the following Lemma:
\begin{lemma}\label{Lem_maxAp_Ap_n=1,k=2}
Let us consider $GT(n,k)$ for the case of $n = 1$ and $k = 2$. Then by Theorem \ref{Thm_minimal}, we can obtain $\\GT(1,2) = \big<\{(2^2 + 1)\cdot 2^1 - (2^2 - 1), (2^2 + 1)\cdot 2^2 - (2^2 - 1), (2^2 + 1)\cdot 2^3 - (2^2 - 1)\}\big>$\\
and $Ap(GT(1,2),7) = \{0,17,34,37,54,71,74\}$ can be verified directly.\\
Hence $\max(Ap(GT(1,2),7)) = 74$ and $F(GT(1,2)) = 74 - 7 = 67$.
\end{lemma}

\newpage
\section{Conclusion}\label{sec_Conclusion}
We summarize all previous results and three main Theorems. For convenience, let
\begin{gather*}
P_i = t_1 + (2^2 - 1)t_2 + \cdots + (2^i - 1)t_i
\end{gather*}
and
\begin{gather*}
Q_i = t_1 s_1 + \cdots + t_i s_i
\end{gather*}
in this section.

First main result is about the minimal system of generators of $GT(n,k)$ and it was already stated in Theorem \ref{Thm_minimal}.

By Lemma \ref{Lem_maxAp_n=0}, Lemma \ref{Lem_maxAp_k=n}, Lemma \ref{Lem_maxAp_2<k<n}, Lemma \ref{Lem_maxAp_k>n} and Lemma \ref{Lem_maxAp_Ap_n=1,k=2}, we obtain $\max(Ap(GT(n,k),s_0))$ and $F(GT(n,k))$ for general case and it is stated in Theorem \ref{Thm_conclusion_1} as the following. Note that the case for $n \neq 0$ and $k = 1$ is in \cite{Rosales2015}.

\begin{theorem}\label{Thm_conclusion_1}
Let $n \in \mathbb{N}$, $k \in \mathbb{N}\texttt{\char`\\}\{0\}$ and let a generalized Thabit numerical semigroup which is associated with $n,k$ be $GT(n,k) = \big<\{(2^k + 1)\cdot 2^{n+i} - (2^k - 1)| i \in \mathbb{N}\}\big>$. Then we can obtain the maximal element of the Ap$\acute{e}$ry set for generalized Thabit numerical semigroups.\\
$\max(Ap(GT(n,k),s_0)) =$
\begin{displaymath}
\left\{ \begin{array}{ll}
s_1 & \text{if} ~~ n = 0,\\
s_n + s_{n+1} & \text{if} ~~ n \neq 0 ~~ and ~~ k = 1,\\
Q_{n-2} + s_{n-1} + s_{n+k} ~~ where ~~ P_{n-2} = 2^{n-1} - 2^k + 2 & \text{if} ~~ 2 \leq k < n,\\
s_1 + s_{2n} & \text{if} ~~ n = k \geq 1,\\
Q_{k-1} + s_k + \cdots + s_{n+k-1} ~~ where ~~ P_{k-1} = 2^n + n & \text{if} ~~ n \neq 0, 2 \neq k > n,\\
74 & \text{if} ~~ n = 1 ~~ and ~~ k = 2.
\end{array} \right.
\end{displaymath}
And by using this result, we can obtain the Frobenius number for generalized Thabit numerical semigroups.\\
$F(GT(n,k)) =$
\begin{displaymath}
\left\{ \begin{array}{ll}
s_1 - 2 & \text{if} ~~ n = 0,\\
s_n + s_{n+1} - s_0 & \text{if} ~~ n \neq 0 ~~ and ~~ k = 1,\\
Q_{n-2} + s_{n-1} + s_{n+k} - s_0 ~~ where ~~ P_{n-2} = 2^{n-1} - 2^k + 2 & \text{if} ~~ 2 \leq k < n,\\
s_1 + s_{2n} - s_0 & \text{if} ~~ n = k \geq 1,\\
Q_{k-1} + s_k + \cdots + s_{n+k-1} - s_0 ~~ where ~~ P_{k-1} = 2^n + n & \text{if} ~~ n \neq 0, 2 \neq k > n,\\
67 & \text{if} ~~ n = 1 ~~ and ~~ k = 2.
\end{array} \right.
\end{displaymath}
\end{theorem}
By Lemma \ref{Lem_A}, Lemma \ref{Lem_maxAp_n=0}, Lemma \ref{Lem_Ap_k<n}, Lemma \ref{Lem_Ap_k=n}, Lemma \ref{Lem_Ap_k>n} and Lemma \ref{Lem_maxAp_Ap_n=1,k=2}, we obtain the Ap$\acute{e}$ry set for $GT(n,k)$ as explicit form and it is stated in Theorem \ref{Thm_conclusion_2} as the following. Note that the case for $n \neq 0$ and $k = 1$ is in \cite{Rosales2015}.
\begin{theorem}\label{Thm_conclusion_2}
Let $n \in \mathbb{N}$, $k \in \mathbb{N}\texttt{\char`\\}\{0\}$ and let a generalized Thabit numerical semigroup which is associated with $n,k$ be $GT(n,k) = \big<\{(2^k + 1)\cdot 2^{n+i} - (2^k - 1)| i \in \mathbb{N}\}\big>$. Also, let $A(n,\delta)$ be the set of the sequences $(t_1, \cdots, t_{n+\delta}) \in \{0,1,2\}^{n+\delta}$ which satisfies the following conditions:
\begin{enumerate}
\item $t_{n+\delta} = 0$ or $1$,
\item If $t_j = 2$, then $t_i = 0$ for all $i < j$.
\end{enumerate}
And let $R(n), R_1(n,\delta),R_2(n,\delta),R_3(n,\delta) \subseteq A(n,\delta)$ be the set of the sequences which are defined by the following additional statements:
\begin{enumerate}
\item $R(n)$ satisfies the following conditions: \begin{enumerate}
\item $t_{n+1} \in \{0,1\}$,
\item If $t_n = 2$, then $t_{n+1} = 0$,
\item If $t_n = t_{n+1} = 1$, $t_i = 0$ for all $1\leq i <n$.
\end{enumerate}
\item If $t_{n+\delta} = 1$, $R_1 (n,\delta)$ satisfies the following conditions: \begin{enumerate}
\item $t_n = \cdots = t_{n+\delta - 1} = 0$,
\item $t_1 + (2^2 - 1)t_2 + \cdots + (2^{n-2} - 1)t_{n-2} \leq 2^{n-1} - 2^k + 2$.
\end{enumerate}
\item If $t_{2n} = 1$, $R_2 (n,\delta)$ satisfies the following conditions:
\begin{enumerate}
\item $t_1 = 1$,
\item $t_2 = \cdots = t_{2n-1} = 0$.
\end{enumerate}
\item If $t_{\delta + 1} = t_{\delta + 2} = \cdots = t_{n+\delta + 1} = 1$, $R_3 (n,\delta)$ satisfies the following inequality:
\begin{gather*}
t_1 + (2^2 - 1)t_2 + \cdots + (2^{\delta} - 1)t_{\delta} \leq 2^n + n.
\end{gather*}
\end{enumerate}
Then, we have the explicit form of the Ap$\acute{e}$ry set for generalized Thabit numerical semigroups.
\begin{displaymath}
Ap(GT(n,k),s_0) = \left\{ \begin{array}{ll}
\{0, s_1\} & \text{if} ~~ n = 0,\\
\{Q_{n+1} | (t_1,\cdots, t_{n+1}) \in R(n)\} & \text{if} ~~ n \neq 0 ~~ and ~~ k = 1,\\
\{Q_{n+\delta} | (t_1,\cdots, t_{n+\delta}) \in R_1(n,\delta)\} & \text{if} ~~ 2 \leq k < n,\\
\{Q_{n+\delta} | (t_1,\cdots, t_{n+\delta}) \in R_2(n,\delta)\} & \text{if} ~~ n = k \geq 1,\\
\{Q_{n+\delta} | (t_1,\cdots, t_{n+\delta}) \in R_3(n,\delta)\} & \text{if} ~~ n \neq 0, 2 \neq k > n,\\
\{0,17,34,37,54,71,74\} & \text{if} ~~ n = 1 ~~ and ~~ k = 2.
\end{array} \right.
\end{displaymath}
Hence, you can get $g(GT(n,k))$ by using the formula $g(S) = \frac{1}{x} (\sum_{w \in Ap(S,x)} w) - \frac{x-1}{2}$ in Lemma \ref{lem_F_g}.
\end{theorem}
Finally, we suggest some open problems related to this paper.
\begin{open problem}
The Frobenius problem for \textit{more generalized form of Thabit numerical semigroups} is still open, such as $ \big< \{(2k + 1)\cdot 2^{n+i} - (2k - 1) | i \in \mathbb{N}\}\big>$ for any $n \in \mathbb{N}$ and $k \in \mathbb{N}\texttt{\char`\\}\{0\}$. We have tried to solve the problem but it needs far more laborious calculation. We hope that there exist some more efficient methods to solve it. 
\end{open problem}
\begin{open problem}
Let $A$ be the minimal system of generators of $S$ and $B$ be a proper subset of $A$. In other words, $\big<A\big> = S$ and $\big<B\big> \subsetneq S$ but $F(B) = F(A)$ may be satisfied in this condition. For example, let $S = \big< \{(2^3 + 1)\cdot 2^{2+i} - (2^3 - 1) | i \in \mathbb{N}\}\big>$. Then the minimal system of generators of $S$ is $A = \big< \{(2^3 + 1)\cdot 2^{2+i} - (2^3 - 1) | i \in \{0,1,2,3,4\} \}\big>$ since $4 = 2 + 3 - 1$. By letting $B = \big< \{(2^3 + 1)\cdot 2^{2+i} - (2^3 - 1) | i \in \{0,1,2,3\} \}\big>$, $F(B) = F(A) = 1095$ can be easily observed. Finding the condition of the proper subset $B$ of the minimal system of generators $A$ such that $F(B) = F(A)$ for general numerical semigroup $A$ is an interesing problem.
\end{open problem}

\section*{ACKNOWLEDGEMENT}
The author would like to thank the editors and referees for the helpful suggestions. The author also thanks to lab colleagues (Youngwoo Kwon and Myunghyun Jung) for their support during the course of this study. Finally, The author would like to thank his advisor, Donggyun Kim, for his guidance on this paper.
\section*{Conflicts of interest}
The author declares that there is \bf{No conflict of interests} regarding the publication of this article.

\end{document}